\documentclass[10pt]{article}

\usepackage{amsmath}
 \usepackage{amsfonts,amssymb,graphicx}
\usepackage[colorlinks=true,linkcolor=blue,citecolor=red]{hyperref}

\newtheorem{prop}{Proposition}
\newtheorem{rema}{Remark}
\newtheorem{defi}{Definition}
\newtheorem{lemm}{Lemma}
\newtheorem{theo}{Theorem}
\newtheorem{coro}{Corollary}

\newcommand{\X}[1][]{\ensuremath{{\mathbb{S}^{#1}} }}

\newcommand{\C}[1][]{\ensuremath{{\mathbb{C}^{#1}} }}
\newcommand{\R}[1][]{\ensuremath{{\mathbb{R}^{#1}} }}

\renewcommand{\S}[1][]{\ensuremath{{\mathbb{S}^{#1}} }}
\renewcommand{\H}[1][]{\ensuremath{{\mathbb{H}^{#1}} }}

\newcommand{\M}{{\cal M}}
\newcommand{\N}{{\cal N}}
\newcommand{\ar}{{\texttt{arctan}\eps}}
\newcommand{\co}{{\texttt{cos}\eps}}
\newcommand{\si}{{\texttt{sin}\eps}}
\renewcommand{\L}{{\cal L}}
\newcommand{\s}{{\cal S}}
\newcommand{\G}{\ensuremath{\mathbb{G}}}
\newcommand{\J}{\ensuremath{\mathbb{J}}}
\newcommand{\<}{\langle}
\renewcommand{\>}{\rangle}
\newcommand{\ga}{\gamma}
\newcommand{\pa}{\partial}
\newcommand{\ka}{\kappa}

\newcommand{\eps}{\epsilon}

\newcommand{\la}{\lambda}
\newcommand{\be}{\beta}
\newcommand{\de}{\delta}

\date{}
\title{Spaces of geodesics of pseudo-Riemannian space forms and normal congruences of hypersurfaces}
\author{ Henri Anciaux\footnote{The author is supported by CNPq (PQ 302584/2007-2) and Fapesp (2010/18752-0)}}

\begin{document}
\maketitle

\centerline{\textbf {Abstract}}

\medskip

{\small We describe natural K\"ahler or para-K\"ahler structures of the spaces  of geodesics of pseudo-Riemannian space forms  and relate the local geometry of hypersurfaces of space forms  to that of their normal congruences, or Gauss maps, which are Lagrangian submanifolds.

The space of geodesics $L^{\pm}(\X^{n+1}_{p,1})$ of a pseudo-Riemannian space form $\X^{n+1}_{p,1}$ of non-vanishing curvature enjoys a K\"ahler or para-K\"ahler structure $(\J,\G)$ which is in addition Einstein. Moreover, in the three-dimensional case, $L^{\pm}(\X^{n+1}_{p,1})$  enjoys another K\"ahler or para-K\"ahler structure $(\J',\G')$ which is scalar flat. 
The normal congruence of a  hypersurface $\s$ of $\X^{n+1}_{p,1}$ is a Lagrangian submanifold $\bar{\s}$ of  $L^{\pm}(\X^{n+1}_{p,1})$, and we relate the local geometries of $\s$ and $\bar{\s}.$ In particular $\bar{\s}$ is totally geodesic if and only if $\s$ has parallel second fundamental form.  
In the three-dimensional case, we prove that $\bar{\s}$ is minimal with respect to the Einstein metric $\G$ (resp.\ with respect to the scalar flat metric $\G'$) if and only if it is the normal congruence of a minimal surface $\s$ (resp.\ of a surface $\s$ with parallel second fundamental form); 
  moreover $\bar{\s}$ is  flat  if and only if $\s$ is Weingarten.
}

\medskip 

\centerline{\small \em AMS 2000 MSC: 53C50, 53C25, 53D12, 53C42 \em }

\medskip

\section*{Introduction}

After the seminal paper of N. Hitchin (\cite{Hi}) describing the natural complex structure of the space of oriented straight lines of Euclidean $3$-space,  several invariant structures on the space of geodesics of certain Riemannian manifolds and their submanifolds have recently been explored by different authors  (see \cite{AGR}, \cite{Ge}, \cite{GK1}, \cite{GG1}, \cite{GG2}, \cite{Ho}, \cite{Ki}, \cite{Sa1}, \cite{Sa2}). In \cite{AGK}, a unified viewpoint has been given
 to this question, classifying all invariant Riemannian, symplectic, complex and para-complex structures that may exist on the space of geodesics
 of a number of spaces: the Euclidean and pseudo-Euclidean spaces, the Riemannian and pseudo-Riemannian space forms and the complex and quaternionic space forms. One of the interesting issues about the spaces of geodesics is
 that the normal congruence (or Gauss map) of a one-parameter family of parallel hypersurfaces in some space is a Lagrangian submanifold   of the corresponding space of geodesics.

The purpose of this paper is twofold: first, to give a more precise picture of the structure of the space of geodesics
of pseudo-Riemannian space forms, and second to study in detail the relationships between the pseudo-Riemannian geometry of a one-parameter family of parallel hypersurfaces and that of its normal congruence.

In particular, we describe the natural K\"ahler or para-K\"ahler structure of the space of geodesics of pseudo-Riemannian space forms of non-vanishing curvature and  prove  that the corresponding metric $\G$ is Einstein (Theorem \ref{amb}).
 The space of geodesics of pseudo-Riemannian three-dimensional space forms, which
is four-dimensional, is specific since \textit{(i)} it is the only dimension for which the space of geodesics of flat pseudo-Euclidean spaces enjoys an invariant metric (see \cite{Sa1}, \cite{AGK}), and \textit{(ii)} in the non-flat case it enjoys another natural complex or para-complex structure, which in turns defines a neutral metric $\G'.$ We prove that $\G'$ is  scalar flat and locally conformally flat (Theorem \ref{ambn=2}).

Next we turn our attention to the relation between one-parameter families of parallel hypersurfaces in pseudo-Riemannian space forms and their normal congruences. We first check that an $n$-dimensional geodesic congruence $\bar{\s}$ is Lagrangian if and only it crosses orthogonally a hypersurface $\s$ (Theorem \ref{lagr}), and therefore all the hypersurfaces $\s_t$ parallel to $\s$ and to its polar.
 Given a one-parameter family of parallel hypersurfaces $(\s_t)$ and its normal congruence $\bar{\s},$ we relate the first and second fundamental forms of $(\s_t)$ to those of $\bar{\s}$ (Theorems \ref{geolagrG} and \ref{geolagrG'}).
 These formulas imply several interesting corollaries: $\bar{\s}$ is totally geodesic (either with respect to $\G$ or $\G'$)
 if and only if the hypersurfaces $\s_t$ have parallel second fundamental form; in the three-dimensional case, $\bar{\s}$ is minimal  with respect to $\G$  if and only if  one of the parallel surfaces $\s_t$ is minimal (Corollary  \ref{coroGmini});  $\bar{\s}$ is minimal  with respect to $\G'$  if and only if the parallel surfaces $\s_t$ are totally geodesic (Corollary \ref{coroG'mini}); the induced metric on $\bar{\s}$ is flat if and only if the surfaces $\s_t$ are Weingarten  (Corollary  \ref{coroG'flat}). We also exhibit three families of Lagrangian surfaces which are marginally trapped with respect to $\G$ or $\G'.$
 (Corollary \ref{coromargtrapped}).

\medskip

The papers is organised as follows: Section 1 provides some useful preliminaries and Section 2 gives the precise statements of results;  Section 3 deals with the geometry of the spaces of geodesics while Section 4 is devoted to normal congruences of hypersurfaces.

The author thanks Nikos Georgiou for interesting observations about the early version of this manuscript.

\section{Preliminaries}

\subsection{Hypersurfaces in pseudo-Riemannian space forms} \label{preli1}
Consider the real space $\R^{n+2}$ and  endowed with the canonical pseudo-Riemannian metric of signature $(p,n+2-p),$ where $ 0 \leq p \leq n+1$:
$$ \<.,.\>_p := -\sum_{i=1}^{p}  dx^2_i+\sum_{i=p+1}^{n+2}  dx^2_i,$$
and the $(n+1)$-dimensional quadric
$$\X^{n+1}_{p,\eps}=\big\{ x \in \R^{n+2} \big| \, \<x,x\>^2_p= \eps \big\},$$
where $\eps=\pm 1.$
The metric induced on $\X^{n+1}_{p,\eps}$ by the canonical inclusion $  \X_{p,\eps}^{n+1} \hookrightarrow (\R^{n+2}, \<.,.\>_p)$ has signature $(p,n+1-p)$ if $\eps=1$ and $(p-1,n+2-p)$ if $\eps=-1,$ and
constant sectional curvature $K=\eps.$ Conversely, it is known (see \cite{Kr}) that any pseudo-Riemannian manifold with constant sectional curvature is, up to a scaling of the metric, locally isometric to one of these quadrics. 
The transformation
$$\begin{array}{cccc} \texttt{A}:&  \R^{n+2} &\to& \R^{n+2} \\ 
&  (x_1,...,x_p,x_{p+1},...,x_{n+2}) & \mapsto &  (x_{p+1},...,x_{n+2},x_{1},...,x_{p+1})
\end{array}$$
defines an anti-isometry of $\X^{n+1}_{p,\eps}$ onto $\X^{n+1}_{n+2-p,-\eps}.$ Is it therefore sufficient to study the case $\eps=1.$
The two Riemannian space forms are  $(i)$ the sphere $ \X^{n+1}:=\X_{0,1}^{n+1,1},$ which is the only compact quadric, and $(ii)$ the hyperbolic space $\H^{n+1}:=\texttt{A} (\X^{n+1}_{n+1,1}) \cap \{ x \in \R^{n+2} |   x_{1} >0\})$ ($\X^{n+1}_{1,-1}$ and $\X^{n+1}_{n+1,1}$ are the only non-connected quadrics).
Analogously, the two Lorentzian space forms are the de Sitter space $d\S^{n+1}:=\X^{n+1}_{1,1}$ and the anti de Sitter space  
$ Ad\S^{n+1}:=\X_{2,-1}^{n+1}=\texttt{A}(\X_{n,1}^{n+1}).$

Let  $\phi : \M^n \to \X^{n+1}_{p,1}$ be a smooth map from an orientable $n$-dimensional manifold $\M^n.$ 
We set $g:=\phi^* \<.,.\>_p$ for the induced metric on $\M^n.$ We shall always assume that $\phi$ is a pseudo-Riemannian immersion, i.e.\ $g$ is non-degenerate. This is equivalent to the existence of a   unit normal vector field along the immersed hypersurface $\s:=\phi(\M^n)$ that we will denote by $N.$
The curvature of $\s$ may be equivalently described by two tensors: the second fundamental form $h$ with respect
to $N$, i.e.\ $h(X,Y) = g( \nabla_X Y , N),$ where $\nabla$ denotes the Levi-Civita connection of $\<.,.\>_{p}$;
 the shape operator defined by $ AX = -dN(X).$ They  are related by the formula: $g(AX,Y) = h(X,Y).$
The shape operator $A$ is not necessarily real diagonalizable since it is symmetric with respect to the possibly indefinite metric $g.$ 
More precisely, $A$ may be of three types: real diagonalizable, complex diagonalizable, or not diagonalizable at all. In the two-dimensional case, 
we shall use the existence of a canonical form for $A,$ i.e. the existence of a frame $(e_1,e_2)$ such that the matrices of  $g$ and  $A$ take a simple form (see \cite{Ma}):

\begin{itemize}
	\item[-]  real diagonalizable case:
	 $$g= \left( \begin{array}{cc} \eps_1 & 0 \\ 
0 & \eps_2
\end{array} \right) \quad \mbox{ and } \quad 
A=\left( \begin{array}{cc}
 \ka_1 & 0 \\ 
0 & \ka_2
\end{array} \right),$$
with $\eps_1,\eps_2 =\pm 1$;

	\item[-] complex diagonalizable case
	 $$g= \left( \begin{array}{cc} -1& 0 \\ 
0 & 1
\end{array} \right) \quad \mbox{ and } \quad 
A=\left( \begin{array}{cc}
 H & \la \\ 
-\la & H
\end{array} \right),$$
with non-vanishing $\la$; 

	\item[-] non diagonalizable case:
	$$g= 
\left( \begin{array}{cc}
 0& 1 \\ 
1 & 0
\end{array} \right) \quad \mbox{ and }\quad 
A=\left( \begin{array}{cc}
 H & 1 \\ 
0 & H
\end{array} \right).$$
\end{itemize}

\subsection{Parallel hypersurfaces }
It will be convenient to introduce some notation: we set  $(\co, \si):=(\cos, \sin)$ if $\eps=1$ and 
 $(\co, \si):=(\cosh, \sinh)$ if $\eps=-1.$
Given $t \in \R,$ the image of
$$\phi_t:= \co (t) \phi + \si (t) N,$$
when an immersion, is parallel to $\s.$
 When $A$ is invertible, the map $N: \M^n \to \S^{n+1}_{p,\eps},$ where $\eps:=|N|^2_p,$ is an immersion and its image $\s':=N(\M^n)$ is called the \em 
 polar \em of $\s.$
If $\eps=1,$ we have $\phi_{\pi/2}=N,$ hence the polar of $\s$ is parallel to $\s.$ If $\eps=-1,$ $\s' \in \X^{n+1}_{p,-1}= \texttt{A}(\X^{n+1}_{n+2-p,1}).$
In all cases, a unit normal vector field  along $\s_t=\phi_t(\M^n)$ is
$$N_t:=\co (t) N -   \eps \si (t) \phi,$$
 which, when an immersion, is parallel to $\s'.$

\begin{lemm} \label{L1}
Let $\phi:\M^2 \to \X^3_{p,1}$ an immersion with mean curvature $H$ and  Gaussian curvature $K$, which satisfies the following linear Weingarten equation
$$ \left| \frac{2H}{K-\eps} \right|= C,$$
where $C \in [0,  \infty] ,$ and $(\eps,C) \neq (-1,1).$
Then there exists a minimal immersed hypersurface which is parallel to $\s:=\phi(\M^2)$ or to its polar $\s'.$
\end{lemm}

\noindent \textit{Proof.}
We first compute
$$d\phi_t =\co(t) d\phi + \si (t) dN = (\co (t) Id - \si(t) A) \circ d\phi.$$
Observe that $\phi_t$ is an immersion if and only if $\co(t)Id - \si(t)A$ is invertible. When this is the case, we have
$$ dN_t= \co(t) dN - \eps \si(t) d\phi = (-\co (t) A -\eps \si(t) Id) \circ d\phi,$$
so
$$ A_t= -dN_t = (\co (t) A +\eps \si(t) Id)\circ (\co (t) Id - \si(t) A)^{-1}.$$ 
Settting $$A:=\left( \begin{array}{cc}
 H+a & c \\ 
b & H-a
\end{array} \right),$$
we obtain
$$ \co (t) Id - \si(t) A=\left( \begin{array}{cc}
 \co(t) - \si(t) H+a & - \si (t)c \\ 
 - \si(t) b & \co (t) - \si (t) H-a
\end{array} \right).$$
Hence $H_t= tr A_t$ vanishes if and only if the following vanishes as well:
\begin{eqnarray*}&&\co (t)(H+a) + \eps \si (t) )(\co (t)- \si (t) (H-a) \\
&&+ \co (t)(H+a) + \eps \si (t) )(\co (t)- \si (t) (H-a)
 - \co (t)b \si (t) c -  \co (t)c \si (t) b\\
 &=& \co^2 (t) (H+a +H-a) - \eps \si^2(t)(H-a+H+a) \\
 &&+ \co (t) \si (t)\left(2\eps -2 (H-a)(H+a)-2bc \right) \\
 &=&\co(2t) 2H + \si(2t) (\eps - K).
\end{eqnarray*}
Thus, if $\eps=1$ we get the vanishing of $H_{t_0}$ setting $t_0:=\frac{1}{2}\tan^{-1} \left(\frac{2H}{K-1} \right).$
If $\eps=-1$ and $\left|\frac{2H}{K+1} \right|<1,$ the same occurs with 
$t_0:=\frac{1}{2} \tanh^{-1} \left(\frac{2H}{K+1} \right).$ Finally, if $\eps=-1$ and $\left|\frac{2H}{K+1} \right|>1,$ 
we easily check that
$N_{t_0}:=\cosh(t_0) N - \eps  \sinh (t_0) \phi,$ where $t_0:=\frac{1}{2} \coth^{-1} \left(\frac{2H}{K+1} \right),$ is minimal. This completes
the proof.

\medskip

We shall denote by $\ar $ the integral of the map $\frac{1}{1 + \eps t^2 }$, i.e.:
$$ \ar(t) = \left\{ \begin{array}{l} \tan^{-1}(t) \mbox { if } \eps=1,\\
\tanh^{-1}(t) \mbox{ if }  \eps=-1, |t| <1, \\ 
 \coth^{-1}(t) \mbox{ if } \eps=-1, |t| >1. \end{array} \right. $$
 The only  property of $\ar$ we will need is the following:  
 $$\ar(a) + \ar(b)= \ar \left(\frac{a+b}{1 - \eps a b} \right).$$

\subsection{Lagrangian submanifolds}
We first recall the definition of a Lagrangian submanifold:
\begin{defi}
Let $({\cal N}, \omega)$ a $2n$-dimensional symplectic manifold. 
An immersion $\phi: \M^n \to {\cal N}$ is said to be \em Lagrangian \em if $\phi^* \omega =0.$
\end{defi}

We refer the reader to  \cite{AMT} or \cite{CFG} for material  about
  para-complex geometry (sometimes referred to as \em split-complex \em or \em bi-Lagrangian \em geometry).
By a \em pseudo-K\"ahler \em or a \em para-K\"ahler \em manifold, we mean  a manifold equipped with a complex or para-complex structure $\J$ and a compatible pseudo-Riemannian metric $\G,$ i.e.\ such that $\G(\J.,\J.)=\eps \G(.,.).$ Here, $\eps=1$ in the complex case and $\eps=-1$ in the para-complex case. In other words $\J$ is an isometry in the complex case and an anti-isometry in the para-complex case. 
It is furthermore required that the  \em symplectic form \em $\omega:=\eps \G(\J.,.)$  be closed\footnote{Of course the factor $\eps$ is unessential here and is put in order to simplify further exposition. In particular, this convention allows to recover, in the case of $\R^2,$ the "natural" objects $\G:= dx^ 2 + \eps dy^2, $ $\J (\pa_x,\pa_y) := (\pa_y,-\eps \pa_x)$ and $\omega:=dx \wedge dy$.}.
 Observe that the metric $\G$ is determined by the pair $(\J,\omega)$ via the equation $\G:=\omega(.,\J.)$.

It is  well known that the extrinsic curvature of a Lagrangian submanifold in a K\"{a}hler manifold $({\cal N},\J,\G)$ is described by the tri-symmetric tensor
$ {h}(X,Y,Z) := \G(D_X Y, \J Z),$ where $D$ denotes the Levi-Civita connection of $\G$
 (see \cite{An}). 
 It turns out that the same fact holds in the para-K\"ahler case:

\begin{lemm} \label{tri} Let $\L$ be a non-degenerate, Lagrangian submanifold of a pseudo-K\"ahler or para-K\"ahler
manifold $(\N,\J,\G,\omega).$ Denote by $D$ the Levi-Civita connection of  $\G$.
Then
the map $$ {h}(X,Y,Z):=\G(D_X Y,\J Z)$$ is tensorial and
tri-symmetric, i.e.\ $${h}(X,Y,Z)={h}(Y,X,Z)={h}(X,Z,Y).$$

\end{lemm}

\noindent \textit{Proof.}
 The tensoriality of $h$ and its symmetry
with respect to the first two slots follow from the tensiorality and the symmetry
of the second fundamental form. It remains to prove for example that
$h(X,Y,Z)=h(X,Z,Y).$ From the Lagrangian assumption we have
$\omega(Y,Z)=\eps \G(\J Y,Z)=0.$ Differentiating in the $X$ direction
gives, using the fact that $\J$ is parallel,
\begin{eqnarray*} 0&=&X (\G(\J Y,Z))=G(D_X \J Y,Z)+ \G(\J Y,D_X Z)\\
&=&\G(\J D_X Y,Z)+ \G(\J Y,D_X Z)\\
&=&-\G(D_X Y,\J Z)+ \G(\J Y, D_X Z )\\
 &=& -\G(D_X Y,\J Z)+ \G(\J Y, D_X Z)\\
 &=& -h(X,Y,Z)+h(X,Z,Y),
 \end{eqnarray*}
and the proof is complete. 

\bigskip


\section{Statement of results}

\subsection{Structures of the space of geodesics of pseudo-Riemannian space-forms}

Let $x$ be a point of $\X^{n+1}_{p,1}$ and $v \in T_x \X^{n+1}_{p,1} = x^{\perp}$ a unit vector tangent to $x.$ 
Setting $\eps:=|v|^2_p,$ 
the unique geodesic $\ga$ of $\X^{n+1}_{p,1}$ passing through $x$ with velocity $v$ is the periodic curve 
parametrized by $\ga(t)=\co (t) x + \si (t) v$. 
 
The set $L^+(\X^{n+1}_{p,1})$ of positive oriented geodesics of $\X^{n+1}_{p,1}$ identifies with the Grassmannian $Gr^+({n+2},2)$
of oriented two-planes of $\R^{n+2}$ with positive induced metric,
while the set $L^-(\X^{n+1}_{p,1})$ of negative  oriented geodesics of $\X^{n+1}_{p,1}$ identifies  with the Grassmannian $Gr^-({n+2},2)$
of oriented two-planes of $\R^{n+2}$ with indefinite induced metric:
$$L^+(\X^{n+1}_{p,1}) \simeq Gr_p^+({n+2},2) \simeq  \left\{ x \wedge y \in \Lambda^2(\R^{n+2}) \big| (x,y) \in T\X^{n+1}_{p,1}, \<y,y\>_p=1 \right\},$$
$$L^-(\X^{n+1}_{p,1}) \simeq Gr_p^-(n+2,2) \simeq  \left\{ x \wedge y \in \Lambda^2(\R^{n+2}) \big| (x,y) \in T\X^{n+1}_{p,1}, \<y,y\>_p=-1 \right\}.$$
Observe that the anti-isometry $\texttt{A}$ induces a canonical one-to-one correspondence between $Gr_p^-(n+2,2)$ and $Gr_{n+2-p}^-(n+2,2),$
hence between $L^-(\X^{n+1}_{p,1})$ and $L^-(\X^{n+1}_{n+2-p,1}) =L^-(\texttt{A}(\X^{n+1}_{p,-1})).$

We may regard $L^+(\X^{n+1}_{p,1})$ and $L^-(\X^{n+1}_{p,1})$ as two submanifolds of the pseudo-Euclidean space 
$$\Lambda^2(\R^{n+2}):= \mbox{Span} \left\{ e_i \wedge e_j , 1 \leq i < j \leq n+2 \right\} \simeq \R^{\frac{(n+2)(n+1)}{2}},$$
 where $(e_1,...,e_{n+2})$ denotes the canonical basis of $\R^{n+2}.$
This viewpoint allows us to define in a natural way several structures on $L^\pm(\X^{n+1}_{p,1})$: first, we use the fact that
$\Lambda^2(\R^{n+2})$ is equipped with the flat pseudo-Riemannian metric
$$\<\< x \wedge y , x' \wedge y'\>\> := \<x,x'\>_p \<y,y'\>_p -  \<x,y'\>_p \<y,x'\>_p ;$$
we shall denote by $\G$ the induced metric on $L^\pm(\X^{n+1}_{p,1})$, i.e. $\G = \iota^*\<\<.,.\>\>,$ where
$ \iota : L^\pm(\X^{n+1}_{p,1}) \to \Lambda^2(\R^{n+2})$ is the canonical inclusion. Second, observe that
a positive (resp.\ indefinite) oriented plane is equipped with a canonical complex (resp.\ para-complex) structure ${\rm J}$.
Explicitly, given $\bar{x}=x \wedge y \in  Gr_p^\pm({n+2},2),$ with $|x|_p^2=1$ and $|y|^2_p=\eps,$ we set
${\rm J}x=y$ and ${\rm J}y = -\eps x.$ In particular, ${\rm J}^2 = \eps {\rm Id}.$ On the other hand, a tangent vector
to $\iota(L^\pm(\X^{n+1}_{p,1}))$ at the point $\bar{x}$ takes the form $x \wedge X + y \wedge Y,$ where $X,Y \in \bar{x}^\perp.$
 We then define:
 $$ \J (x \wedge X + y \wedge Y) := ({\rm J}x) \wedge X + ({\rm J}y) \wedge Y = y \wedge X  - \eps x \wedge Y.$$
It is straightforward that $\left. \J^2 = \eps \rm{Id} \right|_{\bar{x}},$ i.e.\ $\J$ is an almost complex or para-complex structure.

\begin{theo} \label{amb}
$(L^+(\X^{n+1}_{p,1}),\J,\G)$ is a $2n$-dimensional pseudo-K\"ahler manifold with signature $(2p,2(n-p))$ and
 $(L^-(\X^{n+1}_{p,1}),\J,\G)$ is a $2n$-dimensional para-K\"ahler manifold, hence with neutral signature $(n,n)$.
 In both cases, the metric $\G$ is Einstein, with  scalar curvature $\bar{S}=\eps 2n^2,$ and is never conformally flat.
\end{theo}

\begin{rema}\em    It is not difficult to check that $\G$ and $\J$ are invariant under the natural action of the group $SO(n+2-p,p)$ 
of isometries of $\X^{n+1}_{p,1}.$ 
Such invariant structures have been studied with the Lie algebra formalism in \cite{AGK}, where in particular it
is proved  that such an invariant 
  pseudo-Riemannian metric and complex or para-complex structure  are unique on 
$L^\pm(\X^{n+1}_{p,1})$, for $n \geq 3.$
 The fact that $\G$ is Einstein has been proved in \cite{Le} in the spherical case. \em \end{rema}
 
 \begin{rema} \em \label{remahdepend} The complex structure of $L^+(\X^{n+1}_{p,1})$ may be alternatively described by identifying 
 $L^+(\X^{n+1}_{p,1})$ with the hyperquadric 
 $$\Big\{ [z_1 : ... : z_{n+2}] \,  \Big| \, \, -\sum_{i=1}^{p}  z^2_i+\sum_{i=p+1}^{n+2}  z^2_i = 0 \Big\}$$
  of the pseudo-complex projective space $\C\mathbb P^{n+1}_p$ (see \cite{Os}). \em
\end{rema}

In the three-dimensional case, $L^{\pm}(\X^{3}_{p,1})$ enjoys other natural 
 structures, which may de defined as follows: since the orthogonal two-plane
$\bar{x}^\perp$ admits a canonical orientation (that orientation compatible with the orientations of $\bar{x}$
and $\R^4$), it enjoys a canonical complex or para-complex structure ${\rm J}'$ (depending of whether
the induced metric on $\bar{x}^\perp$ is positive or indefinite). Hence we set
$$\J' (x \wedge X + y \wedge Y) := x \wedge ({\rm J}' X) + y \wedge ({\rm J}'Y). $$
We therefore get another almost complex or para-complex structure on $L^{\pm}( \X^{3}_{p,1}).$
Finally, we introduce one more tensor:
we want to define a pseudo-Riemannian structure $\G'$ on $L^{\pm}(\X^{3}_{p,1})$ with the requirement that 
the pair $(\J',\G')$ induces the same symplectic structure, up to sign, than that of $(\J,\G)$. In other words, we require that $\omega(.,.):=\eps' \G'(\J'.,.)$ be the same that   $\omega(.,.):=\eps \G(\J.,.).$
Hence, we must have:
$$\G'= \omega(.,\J'.)=\epsilon \G(\J.,\J'.)=-\eps  \G(.,\J \circ \J'.) .$$
It turns out that this defines another K\"ahler or para-K\"ahler structure:

\begin{theo} \label{ambn=2}
The two-form $\G':=-\eps\G(.,\J' \circ \J.)$ is symmetric and therefore defines a pseudo-Riemannian metric on $L^{\pm}(\S^3_{p,1}).$  The Levi-Civita connection of $\G'$ is the same than that of  $\G$ and the structures $(\J,\G)$ and $(\J',\G')$ share the same symplectic form $\omega.$
Moreover,
$(L^+(\S^{3}),\J',\G')$, $(L^+(\H^{3}),\J',\G')$, $(L^-(d\S^3),\J',\G')$ and $(L^+(Ad\S^3),\J',\G')$ are  pseudo-K\"ahler manifolds 
while $(L^+(d\S^{3}),\J',\G')$ and $(L^-(Ad\S^{3}),\J',\G')$ are para-K\"ahler manifolds.
In all cases,  the metric $\G'$ has neutral signature $(2,2)$, is scalar flat and  locally conformally flat. 
\end{theo}

\begin{rema} \em The properties of $\G'$  have been derived in \cite{GG1} in the case of hyperbolic space.
\em
\end{rema}
The fact that  $(\J,\G)$ and $(\J',\G')$ share both the same Levi-Civita connection and symplectic form implies that they also share some distinguished classes of submanifolds:

\begin{coro} \label{corohdepend}
Lagrangian surfaces,  flat and totally geodesic submanifolds in $L^{\pm}(\X^{3}_{p,1})$ are the same for $(\J,\G)$ and $(\J',\G').$
\end{coro}

\begin{rema} \em
In some cases, these invariant structures may be defined in a more intuitive way. For example, using the direct sum of self-dual and anti self-dual bivectors in $(\Lambda^2((\R^4, \<.,.\>_0), \<\<.,.\>\>) \simeq (\R^6,\<.,.\>_0)$, one can prove that $L^+(\S^3) \simeq \S^2 \times \S^2$ and that
$$ \G = \<.,.\>_0 \oplus  \<.,.\>_0  \quad \quad  \G' = \<.,.\>_0 \oplus - \<.,.\>_0$$
$$ \J = \left( \begin{array}{cc} j & 0 \\ 0 & j \end{array} \right) 
 \quad \quad \J' = \left( \begin{array}{cc} j & 0 \\ 0 & -j \end{array} \right),$$
where $(j,\<.,.\>_0)$ is the canonical K\"ahler structure of $\X^2$ (see \cite{CU}).
Analogously, in $(\Lambda^2(\R^4, \<.,.\>_2), \<\<.,.\>\>) \simeq (\R^6,\<.,.\>_3)$ the Hodge operator is para-complex and we still have a direct sum
of self-dual and anti self-dual bivectors. A computation then shows that $L^+(A d\S^3) \simeq \H^2 \times \H^2$ and  $L^-(Ad\S^3) \simeq d\S^2 \times d\S^2,$ and again we could describe $(\J,\G)$ and $(\J',\G')$ as product structures built from the canonical K\"ahler and para-K\"ahler structures of $\H^2$ and $d\S^2$ respectively.
On the other hand, the Hodge operator being complex in $(\Lambda^2(\R^4, \<.,.\>_1), \<\<.,.\>\>) \simeq (\R^6,\<.,.\>_2)$, there is no natural direct sum of it into eigen-spaces, and it does not seem possible a priori to describe $L^{\pm}(d\X^3)$ and $L^+(\H^3)$  as a Cartesian product of two surfaces. 
\em \end{rema}

\begin{rema} \em
Since the two complex or para-complex structures $\J$ and $\J'$ commute, their composition $\J'' := \J \circ \J'$ defines one more invariant structure:
if $\J$ and $\J'$ are both complex or both para-complex, then $\J''$ is complex, and if $\J$ and $\J'$ are of different types, $\J''$ is para-complex.
The two-form $\G'':= \omega(.,\J'')$ is not symmetric, so there is no pseudo- or para-K\"ahler structure associated to $\J''.$

Observe also that the triple $(\J,\J',\J'')$ is \em not \em a para-quaternionic structure, since $\J$ and $\J'$ commute rather than anti-commute.
The case $L^-(Ad\S^3)$ excepted, this triple is what is called an \em almost product bi-complex \em structure in \cite{Cr}.
 \em

\end{rema}

\subsubsection*{Structures on $L^{\pm}(\X^3_{p,1})$}
  
 \hspace{-5em}  \begin{tabular}{| l || l || c | c | c | c | c |}
     \hline
     Space form & Space of geodesics
     & $(\eps,\eps')$  & Signature of $\G$ & $\J$ & $\J'$  & $\J''$ \\ \hline \hline
   $\X^3_{0,1}=\S^3 $ & $L(\X^3)$ & $(1,1)$  & $(+,+,+,+)$ & complex & complex & para \\ \hline 
   $\X^3_{1,1} = d\X^3$ & $L^+(d\X^3)$ & $(1,-1)$ &  $(+,-,+,-)$ & complex & para & complex \\ \hline 
   & $L^-(d\X^3) \simeq L^-(\H^3) $ & $(-1,-1)$ &  $(+,+,-,-)$ & para  & complex & complex \\ \hline 
  $\X^3_{2,1} \simeq Ad\X^3$ &  $L^+(Ad\X^3)$ & $(1,1)$ &  $(-,-,-,-)$ & complex & complex & para \\ \hline 
   & $L^-(Ad\X^3)$ & $(-1,-1)$ & $(+,-,-,+)$ & para & para & para \\ \hline 
  $\X^3_{3,1} \simeq \H^3$ &  $L^-(\H^3) \simeq L^-(d\S^3)$ & $(-1,1)$ &  $(-,-,+,+)$ & para & complex & complex \\ 
     \hline
   \end{tabular}
\subsection{Normal congruences of immersed hypersurfaces as Lagrangian submanifolds}
\begin{defi}
Let $\s$ be an immersed surface of  pseudo-Riemannian space form $\X_{p,1}^{n+1}$ 
with  unit normal vector
$N.$ The \em normal congruence \em (or \em Gauss map\em) $\bar{\s}$ of  $\s$ is
 set of geodesics crossing $\s$ orthogonally in the direction $N.$
\end{defi}

\begin{theo} \label{lagr} Let $\phi$ be a pseudo-Riemannian immersion of an orientable manifold $\M^n$ in pseudo-Riemannian space form $\X_{p,1}^{n+1}$ 
with  unit normal vector
$N.$ 
Then the normal congruence of $\s:=\phi(\M^n)$ is the image  of the map $\bar{\phi} : \M^n \to L^\pm(\X^{n+1}_{p,1})$ defined by $\bar{\phi} =\phi \wedge N.$ When $\bar{\phi}$ is an immersion, it is Lagrangian with respect to $\omega.$ 
In this case, $\bar{\s}$ is also the normal congruence of the hypersurfaces parallel to $\s$ and to its polar $\s'.$
Conversely, let $\bar{\phi} : \M^n \to L^\pm(\X^{n+1}_{p,1})$  be an immersion of a simply connected
$n$-manifold. Then $\bar{\s}$ is the normal congruence of an immersed hypersurface of $\X^{n+1}_{p,1}$ 
if and only if $\bar{\phi}$ is Lagrangian. 
\end{theo}
In view of this result, it is natural to relate the geometry of a Lagrangian submanifold to that of the corresponding hypersurface of $\X_{p,1}^{n+1}.$

\begin{theo} \label{geolagrG} 
Let $\phi$ be a pseudo-Riemannian immersion of an orientable manifold $\M^n$ in pseudo-Riemannian space form $\X_{p,1}^{n+1}$ 
with  unit normal vector
$N.$ 
Set $|N|_p^2 := \eps $, 
 denote by  $A$ the shape operator of $\phi$ with respect to $N$ and by $\nabla^g$ the Levi-Civita
connection of $g.$
Then the induced metric $\bar{g}:=\bar{\phi}^* \G,$  with $\bar{\phi} =\phi \wedge N$, is given by the following formula
$$ \bar{g}=\eps g +g(A.,A.).$$
In particular, $\bar{g}$ is non-degenerate if and only if $ \eps Id  + A^2$ is invertible.

Moreover, the extrinsic curvatures $h$  of $\s:=\phi(\M^n)$ and of $\bar{h}$ of $\bar{\s}:=\bar{\phi}(\M^n)$ are related by the formula
$$ \bar{h} = \eps \nabla^g h.$$
 In particular the normal congruence $\bar{\s}$ is totally geodesic if and only if $\s$ has parallel second fundamental form.
\end{theo}%

\begin{rema} \emph{The fact that the tensor $ \bar{h}$ of $\bar{\s}$ is tri-symmetric is equivalent to the Codazzi equation for the hypersurface
$\s.$ }
\end{rema}

\begin{coro} \label{coroGminidiago} If the shape operator $A$ of $\s$ is real diagonalizable (this is always the case if $\eps'=1$), the mean curvature vector of  $\bar{\s}$ with respect to $\G$ is 
$$\vec{H}=-\frac{\eps}{n} \J \bar{\nabla} \left(\sum_{i=1}^n \ar (\kappa_i)\right), $$ 
 where $\kappa_1, ..., \kappa_n$ are the principal curvatures  of $\s$ and $\bar{\nabla}$ is the gradient with respect to the induced metric $\bar{g}$. In particular, if $\s$ is isoparametric (i.e.\ its principal curvatures are constant) or austere (i.e.\ the set of its principal curvatures is symmetric with respect to $0$), then its normal congruence $\bar{\s}$ is $\G$-minimal. 
\end{coro}

\begin{coro} \label{coroGmini}
If $n=2,$  the mean curvature vector of  $\bar{\s}$ with respect to $\G$ is 
$$\vec{H}=-\frac{\eps}{2} \J \bar{\nabla} \ar \left(\frac{2H}{1-\eps K}\right), $$
where $H$ and $K$ denote the mean curvature and the Gaussian curvature of $\s$ respectively. In particular, $\bar{\s}$  is $\G$-minimal if and only if it is the normal congruence of a minimal surface. 
\end{coro}

\begin{rema} \emph{Corollaries \ref{coroGminidiago} and \ref{coroGmini} have been proved in  \cite{Pa} in the spherical case. The fact that the mean curvature vector takes the form $\vec{H}=\frac{\eps}{n} \J \bar{\nabla} \be,$ where $\be$ is an $\S^1$-valued map, is due to the fact that the metric $\G$ is Einstein (cf \cite{HR}). The map $\be$ is called the \em Lagrangian angle \em of the submanifold $\bar{\s}$.}
\end{rema}

In the three-dimensional case, it is natural to study the pseudo-Riemannian geometry of Lagrangian surfaces of $L^{\pm}(\S^3_{p,1})$ with respect to the metric $\G'$ described in Theorem \ref{ambn=2}.

\begin{theo} \label{geolagrG'}
Let $\phi$ be a pseudo-Riemannian immersion of an orientable surface $\M^2$ in pseudo-Riemannian space form $\X_{p,1}^3$ 
with shape operator $A$ and unit normal vector
$N.$ 
Then the induced metric $\bar{g}':=\bar{\phi}^* \G',$  with $\bar{\phi} =\phi \wedge N$, is given by the following formula
$$ \bar{g}'=g(., (A{\rm J}'-{\rm J}'A).).$$
Moreover,
\begin{enumerate}
	\item[-] If $A$ is real diagonalizable, the metric $\bar{g}'$ is degenerate at  umbilic points of $\s:=\phi(\M^2)$ and indefinite elsewhere;
	the null directions of $\bar{g}'$ are the principal directions of $\s$;
	\item[-] If $A$ is complex diagonalizable, the metric $\bar{g}'$ is everywhere definite;
	\item[-] If $A$ is not diagonalizable, the metric $\bar{g}'$ is everywhere degenerate.
\end{enumerate}

\noindent  When $\bar{g}'$ is not degenerate, the extrinsic curvatures $h$ and $\bar{h}$ of $\s$ and $\bar{\s}:=\bar{\phi}(\M^2)$ are related by the formula
$$ \bar{h} =  \eps \nabla^g h.$$
 In particular the normal congruence $\bar{\s}$ of $\s$ is totally geodesic if and only if $\s$ has parallel second fundamental form.
\end{theo}

\begin{coro} \label{coroG'mini}  
 $\bar{\s}$ is
 $\G'$-minimal if and only it is totally geodesic, i.e. $\s$ has parallel second fundamental form. 
 In in addition  $A$ is real diagonalizable,
   $\s$ is the set of equidistant points to a geodesic of  $\X^3_{p,1}.$ 
\end{coro}

\begin{coro} \label{coroG'flat}
The induced metric $\bar{g}'$ is flat (and the metric $\bar{g}$ as well by Corollary \ref{corohdepend}) if and only if the surface $\s$ is Weingarten, i.e.\ there exists a functional relation $f(H,K)=0$ satisfied by the mean curvature and the Gaussian curvature of $\s.$ 
\end{coro}

\begin{rema} \em 
Corollary \ref{coroG'mini} and \ref{coroG'flat} have been proved  in the case of hyperbolic space in \cite{Ge} and \cite{GG2} respectively. 
 Corollary \ref{coroG'flat} has been proved in the case of Euclidean space in \cite{GK2}.
\em 
\end{rema}

\begin{coro} \label{coromargtrapped} If the shape operator $A$ of $\s$ is not diagonalizable, then its normal  congruence $\bar{\s}$ is a  $\G$-marginally trapped  surface, i.e.\ the  mean curvature vector of $\bar{\s}$ with respect to $\G$ is null. 
If $\s$ is  a tube (i.e.\ the set of equidistant points to an arbitrary curve of $\X^3_{p,1}$)  or a surface of revolution, then its normal  congruence $\bar{\s}$ is a  $\G'$-marginally trapped  surface. \end{coro}
 

\section{The geometry of the space of geodesics}

\subsection{The Einstein metric $\G$ (Proof of Theorem \ref{amb})} \label{G}

\subsubsection{The second fundamental form of $h^{\iota} $ and the complex structure~$\J$}

\begin{prop}
The complex or para-complex structure $\J$ is integrable.
\end{prop}

\noindent \textit{Proof.} Let $\bar{x}:= x \wedge y \in L^\pm(\X^{n+1}_{p,1})$ 
with  $|x|^2_p=1$ and $|y|^2_p=\eps$
and let $(e_1, ...,e_n)$ be an orthonormal basis of the orthogonal complement
of $x \wedge y.$ We set $\eps_i:=|e_i|^2_p$ and  $\eps_{n+i}:=\eps \eps_i.$
 Then an orthonormal basis $(E_a)_{1 \leq a \leq 2n}$ of $ T_{\bar{x}} L^\pm(\X^{n+1}_{p,1})$, with $\G(E_a,E_a)=\eps_a,$  is given by
$$ E_i:= x \wedge e_i, \quad \mbox{ and }\quad E_{n+i}:= y \wedge e_i.$$
Fix the index $i$ and introduce the curve  
$$\ga_i(t):=x \wedge y_i(t):= x \wedge (\co_{n+i} \,  (t) \,  y+ \si_{n+i} \,  (t) \, e_i).$$
In particular $\ga_i(0)=\bar{x}$ and $\ga_i'(0)=E_i.$
Introduce furthermore the following orthonormal frame $\bar{\cal V}=(\bar{v}_1,...,\bar{v}_{2n})$ along $\ga_i$:
\begin{eqnarray*} &&\bar{v}_j(t) := x \wedge e_j, \quad \quad \bar{v}_{n+j}(t):= y_i(t) \wedge e_j, \quad  \mbox{ if } j \neq i,\\
&&\bar{v}_i(t) := x \wedge y_i'(t) \quad \quad \bar{v}_{n+i}(t):= y_i(t) \wedge y_i'(t). \end{eqnarray*}
Observe that $\bar{v}_a(0)=E_a, \, \forall a, 1 \leq a \leq 2n.$ 
Moreover,
\begin{eqnarray*} &&\bar{v}'_j= 0, \quad \quad \bar{v}'_{n+j}(0)=e_i \wedge e_j, \quad \quad \mbox{ if } j \neq i,\\
&&\bar{v}'_i(0) = x \wedge y_i''(0)=-\eps_{n+i} x \wedge y= -\eps_{n+i} \bar{x}, \quad \quad \bar{v}'_{n+i}(0)= y_i(0) \wedge y_i''(0)=0.
\end{eqnarray*}
Since $ \bar{v}'_{n+j}(0)=e_i \wedge e_j$ and $\bar{v}'_i(0) = -\eps_i \bar{x} $ are normal to $L^{\pm}(\X^{n+1}_{p,1})$, we
deduce that the frame $\bar{\cal V}$ is parallel along $\ga_i.$ On the other hand, 
$$ \J \bar{v}_i= \bar{v}_{n+i} \quad \mbox{ and } \quad \J \bar{v}_{n+i} = -\eps  \bar{v}_{i}.$$
 It follows that $D_{E_i} \J =\J D_{E_i}$
 so $\J$ is parallel, and therefore integrable. 
 
 \bigskip

	We now proceed to compute the second fundamental form of the immersion~$\iota$.

	\begin{prop} \label{2ffi} The second fundamental form of the embedding $\iota: L^{\pm}(\X^{n+1}_{p,1}) \to \Lambda^2(\R^{n+2})$ is given by the formula
$$ h^{\iota}( v \wedge V, w \wedge W ) =- \<v,w\>_p \<V,W\>_p \bar{x} +  \varpi(v,w) V \wedge W,$$  
where $\varpi$ is the symplectic form of the plane $\bar{x}$ defined by $\varpi(.,.) = \eps \<{\rm J}.,.\>_p.$
\end{prop}

\noindent \textit{Proof.} We have 
	\begin{eqnarray*}
	 h^{\iota}(E_i,E_j) & = & \big(D_{\bar{\ga}_i'} \bar{v}_j \big)^\perp  =  \bar{v}'_j(0)  =  - \eps_{n+i} \delta_{ij} \bar{x} \\
	  h^{\iota}(E_i,E_{n+j})& = & \big(D_{\bar{\ga}_i'} \bar{v}_{n+j} \big)^\perp  =  \bar{v}'_{n+j}(0)   =  e_i \wedge e_j.
	  \end{eqnarray*}
An analogous compution, using the curve $\ga_{n+i}(t)=(\co_i \,  (t) \,  x - \si_i \,  (t) \, e_i)\wedge y$, implies that  
$$ h^{\iota}(E_{n+i},E_{n+j}) = -   \eps_i \delta_{ij} \bar{x}.$$
	We deduce that, given $V,W \in \bar{x}^{\perp},$
	\begin{eqnarray*}
	 h^{\iota}(x \wedge V, x \wedge W)&=&-  \eps \<V,W\>_p \bar{x}, \\
	  h^{\iota}(x \wedge V, y \wedge W)&=&V \wedge W,\\
		 h^{\iota}(y \wedge V, y \wedge W)&=&-  \eps \<V,W\>_p \bar{x}.
		 \end{eqnarray*}
The claimed formula follows from the bi-linearity of $h^{\iota}.$

 \subsubsection{The curvature of $\G$}
 We use Gauss equation and Proposition \ref{2ffi} in order  to compute the curvature tensor $\bar{R}$ of $\G$:
for $ 1 \leq a,b,c,d \leq 2n,$ we have
$$ \G(\bar{R}(E_a,E_b)E_c,E_d)=
 \<\< h^{\iota}(E_a,E_c) ,h^{\iota}(E_b,E_d) \>\> 
 - \<\< h^{\iota}(E_a,E_d) ,h^{\iota}(E_b,E_c) \>\> .$$
In particular, we calculate
\begin{eqnarray*} \G(\bar{R}(E_i,E_j)E_k,E_l)&=&
 \<\< h^{\iota}(E_i,E_k) ,h^{\iota}(E_j,E_l) \>\> 
 - \<\< h^{\iota}(E_i,E_l) ,h^{\iota}(E_j,E_k) \>\> \\
&=& \eps_{n+i} \eps_{n+j} \<\<\bar{x},\bar{x} \>\> (\de_{ik} \de_{jl}- \de_{il} \de_{kl})\\
&=& \eps \eps_{i} \eps_{j}  (\de_{ik} \de_{jl}- \de_{il} \de_{kl}).\end{eqnarray*}
 This expression vanishes unless $\{k,l\}=\{i,j\}$ and $i \neq j,$ in which case it becomes
   $$\G(\bar{R}(E_i,E_j)E_i,E_j)=-\G(R(E_i,E_j)E_j,E_i)= \eps \eps_{i} \eps_{j} .$$
A similar computation shows that
$$\G(\bar{R}(E_{n+i},E_{n+j})E_{n+i},E_{n+j})=
-\G(\bar{R}(E_{n+i},E_{n+j})E_{n+j},E_{n+i})= \eps \eps_{n+i} \eps_{n+j}=\eps \eps_{i} \eps_{j}. $$
Moreover,
\begin{eqnarray*}
\G(\bar{R}(E_i,E_{n+j})E_{k},E_{n+l}) &=&
\<\< h^{\iota}(E_i,E_k) ,h^{\iota}(E_{n+j},E_{n+l}) \>\> 
 - \<\< h^{\iota}(E_i,E_{n+l}) ,h^{\iota}(E_{n+j},E_k)\>\>\\
 &=& \eps_{n+i} \eps_j \<\< x,x\>\> \de_{ik} \de_{jl}- \<\< e_i  \wedge e_l , e_k \wedge e_j \>\> \\
 &=& \eps_i \eps_j \de_{ik} \de_{jl} - (\eps_i \eps_j \de_{ik} \de_{jl}-\eps_i \eps_k \de_{ij} \de_{lk})\\
 &=& \eps_i \eps_k \de_{ij} \de_{kl}.
\end{eqnarray*}
This expression vanishes unless $(i,k)=(j,l)$, in which case we get
$$ \G(\bar{R}(E_i,E_{n+i})E_{k},E_{n+k})=\eps_i \eps_k.$$
Using the symmetry of $\G(\bar{R}(.,.).,.),$ we have
\begin{eqnarray*}
\G(\bar{R}(E_{n+i},E_{j})E_{n+k},E_{l}) &=&
\G(\bar{R}(E_{j},E_{n+i})E_{l},E_{n+k})\\
&=&  \eps_i \eps_k \de_{ij} \de_{kl}.
\end{eqnarray*}
Morever, 
\begin{eqnarray*}
\G(\bar{R}(E_i,E_{j})E_{n+k},E_{n+l}) &=&
\<\< h^{\iota}(E_i,E_{n+k}) ,h^{\iota}(E_{j},E_{n+l}) \>\> 
 - \<\< h^{\iota}(E_i,E_{n+l}) ,h^{\iota}(E_{j},E_{n+k})\>\>\\
 &=&  \<\< e_i \wedge e_k, e_j \wedge e_l\>\> - \<\< e_i  \wedge e_l , e_j \wedge e_k \>\> \\
 &=& \eps_i \eps_j (\de_{ik} \de_{jl} -  \de_{il} \de_{jl}).
\end{eqnarray*}
This expression vanishes unless $\{k,l\}=\{i,j\}$ and $i \neq j,$ in which case it becomes
   $$\G(\bar{R}(E_i,E_j)E_{n+i},E_{n+j})=
   -\G(\bar{R}(E_i,E_j)E_{n+j},E_{n+i})=  \eps_{i} \eps_{j}. $$
An easy but tedious calculation shows that $\G(\bar{R}(E_a,E_b)E_c, E_d)$ vanishes if exactly one or three of the indices
$a,b,c$ and $d$ belongs to $\{1,...,n\}.$

\medskip

 It is now easy to calculate the Ricci curvature of $\bar{R}$:
 \begin{eqnarray*}
 \overline{Ric}( E_i,E_j  ) & =  & \sum_{a=1}^{2n}  \G^{aa} 
 \G(\bar{R}(E_i,E_a) E_j ,E_a ) \\
 &=& \sum_{k=1}^{n} \left( \G^{kk} 
 \G(\bar{R}(E_i,E_k) E_j ,E_k ) + \G^{n+k,n+k} 
 \G(\bar{R}(E_i,E_{n+k}) E_j ,E_{n+k} )\right)\\
&=& \sum_{k=1,k\neq i}^n   \eps_k (\de_{ij} \eps \eps_k \eps_i) + \sum_{k=1}^n\eps_{n+k} (\de_{ik} \de_{jk}\eps_{i} \eps_{k}) \\
&=& \delta_{ij}( (n-1) \eps \eps_i + \eps \eps_i)\\
&=& \eps n  \, \G_{ij} 
 \end{eqnarray*}
 and 
 \begin{eqnarray*}
 \overline{Ric}( E_{n+i},E_{n+j}  ) & =  &  \sum_{a=1}^{2n}  \G^{aa} 
 \G(\bar{R}(E_{n+i},E_a) E_{n+j} ,E_a ) \\
&=& \sum_{k=1}^{n} \left( \G^{kk} 
 \G(\bar{R}(E_{n+i},E_k) E_{n+j} ,E_k ) + \G^{n+k,n+k} 
 \G(\bar{R}(E_{n+i},E_{n+k}) E_{n+j} ,E_{n+k} )\right)\\
&=& \sum_{k=1}^n  \eps_k (\delta_{ik} \delta_{jk} \eps_i \eps_k ) +\sum_{k=1,k \neq i}^n  \eps_{n+k}  (\delta_{ij} \eps \eps_{n+k} \eps_{n+j} )  \\
&=&\delta_{ij} (\eps_i + (n-1) \eps \eps_{n+i}) \\
&=& \eps n \,  \G_{n+i,n+j} .
 \end{eqnarray*}
 An analogous calculation shows that  $\overline{Ric}( E_i,E_{n+j})$ vanishes. Hence
 the metric $\G$ is Einstein, with constant scalar curvature $\bar{S}=\eps 2n^2.$
 
 \medskip
 
 Finally, since $\G$ is Einstein, the Weyl tensor is given by the formula
\begin{eqnarray*} W^{\G} &=& \G(\bar{R}.,.) - \frac{\bar{S}}{4n(2n-1)} {\G} \circ \G \\
&=& \G(\bar{R}.,.) - \frac{\eps n}{2(2n-1)} {\G} \circ \G.\end{eqnarray*}
It is easily seen, for example, that
${\G} \circ \G(E_i,E_j,E_{n+i},E_{n+j})$ vanishes. On the other hand, if $i \neq j,$  $\G (\bar{R}(E_i,E_j)E_{n+i},E_{n+j})=\eps_i \eps_j,$ so
$W^\G$ does not vanish and therefore $\G$ is never conformally flat.

\subsection{ The scalar-flat metric $\G'$ in dimension $n=2$ (Proof of Theorem \ref{ambn=2})} \label{G'}
We are going the express all the relevant tensors in  the orthonormal basis $(E_1,E_2,E_3,E_4)$ of $T_{\bar{x}} L^{\pm}(\X^{n+1}_{p,1}).$ Observe first that  matrix of 
$\G$ in this basis is $diag( \eps_1, \eps_2,  \eps_3 ,  \eps_4)=diag( \eps_1, \eps_2, \eps \eps_1 , \eps \eps_2)$ and that
$$ \J =\left( \begin{array}{cccc} 0 & 0 & -\eps & 0 \\ 0 & 0 & 0 &-\eps \\ 1 &0 &0&0 \\ 0 & 1 & 0 &0 \end{array} \right),$$
$$ \J' =\left( \begin{array}{cccc} 0 & -\eps' & 0 & 0 \\ 1  &  0 & 0 &0 \\ 0 & 0 &0 & - \eps' \\ 0&0&1&0\end{array} \right).$$
It follows that
$$\eps  \J \J' = \eps  \J' \J =\left( \begin{array}{cccc} 0 & 0 & 0 &  \eps' \\ 0 & 0 & -1  & 0 \\ 0 & -\eps \eps' &0&0 \\  \eps & 0 & 0 &0 \end{array} \right).$$
Hence, taking into account that $\eps'=\eps_1 \eps_2,$ the matrix of the bilinear form $\G':=-\eps  \G(.,\J \circ \J'.)$ in the basis 
$(E_a)_{ 1 \leq a \leq 4}$ is:
$$ \G' =\left( \begin{array}{cccc} 0 & 0 & 0 & \eps_2 \\ 0 & 0 & -\eps_2 & 0 \\ 0 & -\eps_2 &0&0 \\ \eps_2 & 0 & 0 &0 \end{array} \right).$$

The fact that $\G$ and $\G'$ have the same Levi-Civita connection follows from the next lemma:

\begin{lemm} Let $(\N,\G)$ be a pseudo-Riemannian manifold with Levi-Civita connection $D$ and $T$ a symmetric, $D$-parallel $(1,1)$ tensor.
Then the Levi-Civita connection of the pseudo-Riemannian metric $\G'(.,.):=\G(.,T.)$ is $D.$
\end{lemm}

\noindent \textit{Proof.} Elementary using local coordinates and the explicit formula for the Christoffel symbols.

\bigskip

Since
$\G$ and $\G'$ have the same Levi-Civita connection, they have the same curvature tensor $\bar{R}.$ Therefore,
$$\G'(\bar{R}(.,.).,.):=-\eps  \G(\bar{R}(.,.).,\J \circ \J'.).$$
We deduce
\begin{eqnarray*}
-\overline{Ric}'(X,Y)&= &- \sum_{c,d=1}^4 (\G')^{cd}  \G'(\bar{R}(X,E_c) Y, E_d)\\
&=& (\G')^{14} \G(\bar{R}(X,E_1) Y,\eps  \J \circ \J' E_4)
+ (\G')^{23} \G(\bar{R}(X,E_2) Y,\eps  \J \circ \J' E_3) \\
        && + 
          (\G')^{32} \G(\bar{R}(X,E_3) Y,\eps  \J \circ \J' E_2) 
        +   (\G')^{41} \G(\bar{R}(X,E_4) Y,\eps  \J \circ \J' E_1)\\
 &=&\eps_2 \G(\bar{R}(X,E_1) Y,\eps' E_1)  -\eps_2 \G(\bar{R}(X,E_2) Y,- E_2) \\
 &&        -\eps_2 \G(\bar{R}(X,E_3) Y,-\eps \eps'E_3) 
         + \eps_2 \G(\bar{R}(X,E_4) Y, \eps E_4)\\     
&=&  \eps_1 \G(\bar{R}(X,E_1) Y, E_1)  +\eps_2 \G(\bar{R}(X,E_2) Y, E_2) \\
  &&      +\eps_3 \G(\bar{R}(X,E_3) Y, E_3) 
         + \eps_4 \G(\bar{R}(X,E_4) Y, E_4)\\
         &=& \overline{Ric}(X,Y)\\
         &=& \eps 2 \G(X,Y).
\end{eqnarray*}
It follows that the scalar curvature of $\G'$ vanishes: 
\begin{eqnarray*}
\bar{S}' &=& \sum_{a,b=1}^4 (\G')^{ab}  \overline{Ric}'(E_a,E_b)\\
&=&  -2\eps  \sum_{a,b=1}^4 (\G')^{ab} \G_{ab} \\
&=&  0.
\end{eqnarray*}
It may be interesting to point out that the Ricci curvature of $\G'$ is non-negative in the case of
$L(\S^3)$, non-positive in the case of $L^+(Ad\S^3)$, and indefinite in the other cases.

\bigskip

 Finally, since $\G'$ is scalar flat, its Weyl tensor is given by the formula
\begin{eqnarray*} W^{\G'} &=& \G'(\bar{R}.,.) - \frac{1}{2} \overline{Ric}' \circ \G' \\
&=& \G(\bar{R}.,\eps \J \circ \J'.) - \eps \G \circ \G'.
\end{eqnarray*}
 We may calculate, for example, that
\begin{eqnarray*} W^{\G'}(E_1,E_2,E_2,E_4) &=&  \G(\bar{R}(E_1,E_2)E_1,\eps \J \circ \J' E_4)
 - \eps \, \G \circ \G' (E_1,E_2,E_2,E_4)\\
 &=& -\eps' \eps \eps_1 \eps_2 + \eps \, \G(E_2,E_2) \G'(E_1,E_4)\\
 &=& -\eps +\eps \eps_2 \eps_2 \\
 &=&0.
\end{eqnarray*}
It is easily checked in the same manner that the other components of the Weyl tensor vanish. The metric $\G'$ is therefore locally conformally flat.

\section{Normal congruences of hypersurfaces and Lagrangian submanifolds}

\subsection{Lagrangian submanifolds are normal congruences (proof of Theorem \ref{lagr})}
Let  $\phi : \M^n \to \X^{n+1}_{p,1}$ an immersed, orientable hypersurface with non-degenerate metric and unit normal vector $N$
and introduce the map
$$ \begin{array}{lccc}  \bar{\phi} :
 & \M^n 
 &\to&  L^\pm(\X^{n+1}_{p,1})  \\
& x & \mapsto &  \phi(x) \wedge N(x).
\end{array}$$
In the following, we shall  often allow the abuse of notation of identifying a tangent vector $X$ to $\M^n$ with its image $d\phi(X)$, a vector
tangent to $\X^{n+1}_{p,1}$, therefore an element of $\R^{n+2}.$ We furthermore set $\bar{X}:=d\bar{\phi}(X),$ so that 
$$\bar{X}:=d\bar{\phi}({X})=d(\phi \wedge N)(X)=d\phi(X) \wedge N + \phi \wedge dN(X)=X  \wedge N + AX \wedge \phi.$$
It follows that
\begin{eqnarray*} \omega(\bar{X},\bar{Y})&=&\eps \G(\J \bar{X},\bar{Y})\\
&=&\eps \G(X \wedge ({\rm J}N) +  AX\wedge ({\rm J} \phi),Y \wedge N + AY \wedge \phi)\\
& =&\eps \left( \<X,Y\>_p \<{\rm J} N,N\>_p+\<X,AY\>_p \<{\rm J} N,\phi\>_p
 + \<AX,Y\>_p\<{\rm J} \phi,N\>_p +\<AX,AY\>_p\<{\rm J} \phi,\phi\>_p \right) \\
&= &-\<X,AY\>_p +  \<AX,Y\>_p=0, \end{eqnarray*}
so $\bar{\phi}$ is Lagrangian.

\medskip

hence parallel hypersurfaces have the same normal congruence.
%

Conversely, 
let $\bar{\s}$ an $n$-dimensional geodesic congruence, i.e.\ the image of an immersion $\bar{\phi} : \M^n \to L^{\pm}(\X^{n+1}_{p,1})$. We shall investigate under which condition
there exists an hypersurface $\s$ of $\X^{n+1}_{p,1}$ which intersects orthogonally the geodesics $\bar{\phi}(x),$ $\forall x \in \M^n.$
For this purpose  set $\bar{\phi}(x):=e_1(x) \wedge e_2(x)$ with  $|e_1|^2_p=1$ and $|e_2|^2_p=\eps.$ 
Let $\phi: \M^n \to \X^{n+1}_{p,1}$ such that $\phi(x) \in \bar{\phi}(x), \,  \forall x \in \M^n.$
Therefore there exits $t : \M^n \to \S^1,$ such that 
$\phi(x) = e_1(x) \co (t(x)) + e_2(x) \si (t(x)).$ 
Remember that ${\rm J}$ denotes the complex or para-complex structure on $\bar{\phi}(x),$ in particular
${\rm J} \phi= e_2 \co \, (t) - \eps e_1 \si \, (t).$
 
A computation gives
$$ d\phi=  ({\rm J} \phi) dt + de_1 \wedge e_2 + e_1 \wedge de_2.$$
Since $|\phi(x)|_p^2 =1,$ we always have $\<d\phi ,\phi\>_p =0.$ Hence $\s$ intersects
 the geodesic $\bar{\phi}(x)=e_1(x) \wedge e_2(x) = \phi(x) \wedge {\rm J} \phi(x)$ orthogonally  at the point $\phi(x)$ if and only if the following vanishes:
$$\<d\phi, {\rm J} \phi\>_p =  |{\rm J} \phi|_p^2 d t + \<de_1, e_2\>_p \co^2 (t) - \eps \<de_2, e_1\>_p \si^2 (t)
= \eps dt + \<de_1, e_2\>_p.$$
Hence, $\bar{\s}$ is the normal congruence of $\s$ if and  only there exists $t : \M^n \to \S^1$ such that $ \<de_1, e_2\>_p= - \eps dt.$
Since $\M^n$ is simply connected, it is sufficient to have $d \<de_1 , e_2\>_p =0.$
Observe that 
$$d \<de_1 , e_2\>_p(X,Y)=\<de_1(X) , de_2(Y)\>_p-\<de_1(Y) , de_2(X)\>_p.$$
On the other hand, $ d\bar{\phi} = de_1 \wedge e_2 + e_1 \wedge de_2, $ and 
$$\J(de_1 \wedge e_2 + e_1 \wedge de_2)=-de_1 \wedge e_1 + e_2 \wedge de_2,$$
so that
$$ \omega ( d\bar{\phi}(X), d\bar{\phi}(Y))
=\<\<de_1(X) \wedge e_2 + e_1 \wedge de_2(X), -de_1(Y) \wedge e_1 + e_2 \wedge de_2(Y)\>\>$$
$$=-\<de_1(X) , de_2(Y)\>_p+\<de_1(Y) , de_2(X)\>_p.$$
We conclude that $t$, and thus $\phi$ as well, exists if and only if $\bar{\phi}$ is Lagrangian. Of course, the choice of different constants of integration when solving $t$ corresponds to different, parallel hypersurfaces.

\subsection{Geometry of Lagrangian submanifolds with respect to the Einstein metric $\G$ }

 \subsubsection{The induced metric $\bar{g}= \bar{\phi}^* \G$ 
 and the second fundamental form $\bar{h}$ (proof of Theorem \ref{geolagrG})} \label{2ffG}
Using the description of the metric $\G$ given in Section \ref{G}, we have:
\begin{eqnarray*} \bar{g}(X,Y)&=&
\G(\bar{X},\bar{Y}) \\
 &=&\G( X\wedge N + AX \wedge \phi, Y\wedge N + AY \wedge \phi)\\
& =& \<X,Y\>_p \<N,N\>_p - \<X,N\>_p\<Y,N\>_p 
+ \<\phi, Y\>_p\<AX,N\>_p - \<\phi, N\>_p\<AX, Y\>_p \\
&&+ \<X,\phi\>_p\<N,AY\>_p- \<X,AY\>\<N,\phi\>_p+ \<AX,AY\>_p\<\phi,\phi\>_p-\< AY,\phi\>_p\<AX,\phi \>_p \\
& =&  \eps g(X,Y)+ g(AX,AY).\end{eqnarray*}
We now discuss the degeneracy of $\bar{g}:$
suppose there exist $X$ such that 
$$\bar{g}(X,Y)=\eps g(X,Y)+ g(AX,AY)=g(\eps X+A^2 X,Y) $$
 vanishes $\forall \, Y \in T\M.$
Since the metric $g$ is  non-degenerate, it follows that
$\eps X + A^2 X$ vanishes. Hence $\eps Id + A^2$ is not invertible. If $A$ is diagonalizable, the eigenvalues of $A^2$ are non-negative, 
so we must have $\eps=-1.$

Next, denoting by $\nabla$ (resp. $D$) the flat connection of $\R^{n+2}$ (resp. $\Lambda^2(\R^{n+2})$),
 we have
$$D_{\bar{X}} \bar{Y} = (\nabla_X Y) \wedge N +  (\nabla_X AY) \wedge \phi,$$ 
so
\begin{eqnarray*} \bar{h}(X,Y,Z)&=& \G(D_{\bar{X}} \bar{Y}, \J \bar{Z} ) \\
& =& \G(\nabla_X Y \wedge N +  \nabla_X AY \wedge \phi, Z \wedge (JN) +  AZ \wedge (J\phi))\\
&=&  \G(\nabla_X Y \wedge N +  \nabla_X AY \wedge \phi, -\eps Z \wedge \phi +  AZ \wedge N) \\
&=&  \<\nabla_X Y , AZ \>_p \<N,N\>_p - \eps \< \nabla_X A Y, Z\>_p \<\phi,\phi\>_p \\
&=&\eps \Big( h(\nabla_X Y,Z) - \big( X(\<AY,Z\>_p) - \<AY, \nabla_X Z\>_p \big) \Big) \\
&=&\eps \Big( h( \nabla_X Y,Z)  - X (h(Y,Z)) + h(Y,\nabla_X Z) \Big) \\
& =&\eps(\nabla_X h)(Y,Z) .\end{eqnarray*}

\subsubsection{The mean curvature vector in the diagonalizable case (proof of Corollary \ref{coroGminidiago})} \label{HG}

Assume that $A$ is real diagonalizable and let $(e_1,...,e_n)$ be an orthonormal frame $(e_1,...,e_n)$ on $(T \M,g),$ with $\eps_i:={g}(e_i,e_i)$ 
and such that
$Ae_i = \ka_i e_i,$ where $\kappa_1,...,\kappa_n$ are the principal curvatures of $\s.$

We introduce the notation $\omega^i_{jk}:={g}(\nabla_{e_i} e_j , e_k).$ In particular  $\omega^i_{jk}$ is antisymmetric in its lower indices.
It follows  that
$$ \bar{g}(e_i,e_j)=0 \mbox{ if } i \neq j, \quad \mbox{and}\quad \bar{g}(e_i,e_i)= \eps \eps_i + \eps_i \kappa_i^2=\eps_i(\eps+\ka_i^2).$$
Moreover, if $j \neq k,$
$$ \bar{h}(e_i,e_j,e_k)=\eps \Big( h(\nabla_{e_i} e_j,e_k)+h(e_j, \nabla_{e_i} e_k)-e_i(h(e_j,e_k) ) \Big)
= \eps( \ka_k -  \ka_j ) \omega^i_{jk},$$
and
$$ \bar{h}(e_i,e_j,e_j)=\eps \Big( 2h(\nabla_{e_i} e_j,e_j)- e_i(h(e_j,e_j) )\Big)= - \eps \eps_j e_i(\kappa_j).$$
For further use, observe that the tri-symmetry of $\bar{h}$, or equivalently the Codazzi equation of the immersion $\phi$ implies
$$ ( \ka_j -  \ka_i ) \omega^i_{ij} = \eps_j e_i(\kappa_j). $$
Since the basis $(e_1, ...,e_n)$ is orthogonal with respect to the metric $\bar{g}$, we have
\begin{eqnarray*} \G (n\vec{H},\J d\bar{\phi} (e_i) )&=& \sum_{j=1}^n \frac{\bar{h}(e_j,e_j,e_i)}{\bar{g}(e_j,e_j)}\\
&=&  -\sum_{j=1}^n \frac{\eps  \eps_j e_i(\kappa_j)}{\eps_j(\eps+  \ka_j^2)}\\
&=&  -\sum_{j=1}^n \frac{e_i(\kappa_j)}{1+  \eps \ka_j^2 }\\
& =&-\sum_{j=1}^n e_i(\ar (\kappa_j))\\
&=&  e_i (\beta),\end{eqnarray*}
where $\be:=- \sum_{j=1}^n \ar (\kappa_j),$
which implies that $\vec{H}=\frac{\eps}{n} \J \bar{\nabla} \be.$

Clearly the immersion $\bar{\phi}$ is $\G$-minimal if and only the  map $\be$ is constant. This happens of course if the principal curvatures of $\s$ are constant, i.e.\ it is isoparametric. Moreover, if $\s$ is austere, i.e.\ the set of the principal curvatures is symmetric with respect to $0,$ the Lagrangian angle $\be$ vanishes because the function $\ar$ is odd. This completes the proof of Corollary \ref{coroGminidiago}.

\subsubsection{The mean curvature vector in the two-dimensional case (proof of Corollary \ref{coroGmini})} \label{HG2}
Here and in the next section, we shall make use of canonical form of $A$ (see Section \ref{preli1} and \cite{Ma}).

\medskip

\noindent \textbf{The real diagonalizable case}

\noindent We use the computation of  the previous section:
\begin{eqnarray*}\be &=&- \left(   \ar (\kappa_1) + \ar(\kappa_2) \right)\\
       &=&- \ar \left(\frac{\ka_1 + \ka_2}{1 - \eps \ka_1 \ka_1} \right) \\
       &=& - \ar \left(\frac{2H}{1 - \eps K} \right),
       \end{eqnarray*}
which is the required expression of the Lagrangian angle $\be.$ We now prove that if $\be$ is constant, the assumptions of Lemma \ref{L1} are satisfied.
Assume by contradiction that $(\eps, \frac{2H}{K-\eps})=(-1,\pm 1).$ It follows that
$\frac{\ka_1 + \ka_2}{\ka_1  \ka_2 +1}= \pm 1,$ which in turn implies that $|\ka_1|$ or $|\ka_2|=1.$ Therefore, $-Id + A^2$ is not invertible, and the metric $\bar{g}$ is degenerate by Theorem \ref{geolagrG}. Since this situation is excluded a priori, we may use Lemma \ref{L1} and conclude that there exists
a minimal hypersurface parallel to $\s$ or its polar $\s'$, and therefore whose normal congruence is $\bar{\s}.$
  
 \bigskip
 
\noindent \textbf{The complex diagonalizable case} \label{Hcompldiag}

\noindent Recall that there exists an orthonormal frame $(e_1,e_2)$ such that $g(e_1,e_1)=-1$ and $g(e_2,e_2)=1$ and such that the shape operator takes the form
$$A=\left( \begin{array}{cc}
 H & \la \\ 
-\la & H \end{array} \right)$$
with non-vanishing $\la.$ 
 A quick computation shows that
 $$h=\left( \begin{array}{cc}
 -H & -\la \\ 
-\la & H \end{array} \right) \quad \mbox{ and } \bar{g}= 
\left( \begin{array}{cc}
 -\eps - H^2 + \la^2 & -2H \la \\ 
 -2H\la & \eps + H^2 - \la^2
\end{array} \right). $$
Hence, using the fact that
$$ \nabla_{e_1} e_1 = \omega_{12}^1 e_2, \quad \quad \nabla_{e_1} e_2 = \omega_{12}^1 e_1,$$
$$ \nabla_{e_2} e_1 = \omega_{12}^2 e_2, \quad \quad \nabla_{e_2} e_2 = \omega_{12}^2 e_1,$$
 we calculate
  \begin{eqnarray*} \bar{h}_{111}&=&\eps(-2\la \omega_{12}^1  + e_1(H)) \\
 \bar{h}_{112}&=&\eps(- 2\la \omega_{12}^2  + e_2(H))= \eps e_1(\la)\\
 \bar{h}_{122}&=&  -\eps(2 \la \omega_{12}^1  + e_1(H))=\eps e_2(\la)\\
 \bar{h}_{222}&=&-\eps(2\la \omega_{12}^2  + e_2(H)).
 \end{eqnarray*} 
 (The fact that we obtained two different expressions for $\bar{h}_{112} $ and $ \bar{h}_{122} $ accounts for the Codazzi equation). Hence
\begin{eqnarray*} \G (2\vec{H},\J d\bar{\phi}(e_1) )&=& \frac{ (\eps+ H^2-\la^2)(\bar{h}_{111}-\bar{h}_{122})+4H \la \bar{h}_{112}}
{- (\eps +H^2 - \la^2)^2 - 4 H^2 \la^2}\\
&=&
-\frac{2 \eps (\eps+ H^2-\la^2)e_1(H)+ \eps 4H\la e_1(\la)}{1+H^4+\la^4+ 2\eps H^2 -2 \eps \la^2 - 2H^2 \la^2 +4H^2 \la^2 }\\
&=& -\frac{2  (1+ \eps H^2- \eps \la^2)e_1(H)+ \eps 4H\la e_1(\la)}{1+H^4+\la^4 + \eps 2(H^2 - \la^2) + 2 H^2 \la^2 }.
\end{eqnarray*}
In the same way, we get
\begin{eqnarray*}  \G (2\vec{H},\J d\bar{\phi}(e_2) )&=&  -\frac{2  (1+ \eps H^2- \eps \la^2)e_2(H)+ \eps 4H\la e_2(\la)}
{1+H^4+\la^4 + \eps 2(H^2 - \la^2) + 2 H^2 \la^2}.
\end{eqnarray*}
On the other hand, using the fact that  $K=H^2+\la^2,$ 
\begin{eqnarray*}
d\be& =& d\ar \left(\frac{2H}{1-\eps H^2 - \eps \la^2} \right)\\
&=& \frac{2dH(1-\eps H^2 - \eps \la^2) -2H d(1-\eps H^2 - \eps \la^2 )}{(1-\eps H^2 -\eps \la^2)^2 + \eps 4H^2}\\
 &=& \frac{2(1-\eps 2H^2 -\eps 2\la^2  +\eps 4 H^2 )dH+ \eps 4 H\la d\la}
 {1+H^4+\la^2 - \eps 2H^2 -\eps 2 \la^2 + 2 H^2 \la^2+ \eps 4H^2}\\
&=& \frac{ 2(1 + \eps (H^2 -\la^2)dH  +  \eps 4 H\la d\la}{1+H^4+\la^4 + \eps 2(H^2 - \la^2) + 2 H^2 \la^2 }.
\end{eqnarray*}
It follows that $\G(2\vec{H},\J.)= d\be,$ which is equivalent to $2 \vec{H} =\eps \J \bar{\nabla} \be,$ the required formula.
If  $\eps=-1$ we have, using the fact that $\la \neq 0,$
$$\left| \frac{2H}{K+1} \right|=\frac{2|H|}{H^2+\la^2+1} < \frac{2|H|}{H^2+1} \leq 1.$$
Therefore, if $\bar{\s}$ is $\G$-minimal, i.e. $\be$ is constant, we may use again Lemma \ref{L1} to conclude that there exists a minimal surface
parallel to $\phi$ or $N.$ Hence we have proved Corollary \ref{coroGmini} in this complex diagonalizable case.

\bigskip

\noindent \textbf{The non diagonalizable case} \label{nondiagoG}

\noindent Here there  exists a local frame $(e_1,e_2)$ on $(\M^2,g)$ such that
$$g= 
\left( \begin{array}{cc}
 0& 1 \\ 
1 & 0
\end{array} \right) \quad \mbox{ and }\quad 
A=\left( \begin{array}{cc}
 H & 1 \\ 
0 & H
\end{array} \right).$$
 It follows that
 $$ h=\left( \begin{array}{cc}
 0 & H \\ 
H & 1
\end{array} \right)  \quad \mbox{ and }\quad \bar{g}= 
\left( \begin{array}{cc}
 0 & \eps + H^2 \\ 
 \eps + H^2 & 2H 
\end{array} \right). $$
Observe that $\det ( \eps Id + A^2)=(\eps + H^2)^2,$ hence if $\eps=-1$ and $|H| =1,$ the induced metric $\bar{g}$ is degenerate. Therefore we may assume from now on that $(\eps,|H|) \neq (-1,1).$ Hence, using the fact that
$$ \nabla_{e_1} e_1 = \omega_{12}^1 e_1, \quad \quad \nabla_{e_1} e_2 =- \omega_{12}^1 e_2,$$
$$ \nabla_{e_2} e_1 = \omega_{12}^2 e_1, \quad \quad \nabla_{e_2} e_2 =- \omega_{12}^2 e_2,$$
 we calculate
  \begin{eqnarray*} \bar{h}_{111}&=&0 \\
 \bar{h}_{112}&=& e_1(H)= 0 \\
 \bar{h}_{122}&=& -\eps e_2(H)\\
 \bar{h}_{222}&=&-\eps 2 \omega_{12}^2.
 \end{eqnarray*}
 It follows that
\begin{eqnarray*} \G (2\vec{H},\J d\bar{\phi}(e_1) )&=& \frac{ 2 H \bar{h}_{111}-2(\eps+H^2) \bar{h}_{112}}
{ -(\eps+H^2)^2}\\
&=&0
\end{eqnarray*}
and
\begin{eqnarray*}
 \G (2\vec{H},\J d\bar{\phi}(e_2) )&=& \frac{ 2 H \bar{h}_{112}-2(\eps+H^2) \bar{h}_{122} }{- (\eps+H^2)^2}\\
&=&-2\frac{(\eps+H^2)e_2(H)}{(\eps+H^2)^2}\\
&=&2 \frac{e_2(H)}{1+ \eps H^2}
\end{eqnarray*}
 On the other hand, using the fact that  $K=H^2,$ 
\begin{eqnarray*}
d\be& =& d\ar \left(\frac{2H}{1-\eps K} \right)\\
&=& \frac{2dH(1-\eps H^2) -2H d(1-\eps H^2 )}{(1-\eps H^2)^2 + \eps 4H^2}\\
&=&\frac{2(1+\eps H^2)dH}{(1+\eps H^2)^2}\\
&=&\frac{2 dH}{1+\eps H^2}.
\end{eqnarray*}
taking into account that $e_1(H)$ vanishes, we deduce that $\G(2\vec{H},\J.)= d\be,$ which is equivalent to the required formula. If  $\eps=-1$ we have, using the fact that $|H| \neq 1,$
$$\left| \frac{2H}{K+1} \right|=\frac{2|H|}{1+H^2} < 1.$$
Therefore, we may  use Lemma \ref{L1} again and complete  the proof of Corollary \ref{coroGmini}.

\subsection{Geometry of Lagrangian surfaces with respect to the scalar flat metric $\G'$}

\subsubsection{The metric $\bar{g}'$ and the second fundamental form $\bar{h}$ (proof of Theorem \ref{geolagrG'})}
Using the description of the metric $\G'$ given in Section \ref{G'}, we have:
\begin{eqnarray*} \bar{g}'(X,Y)&=&
\G'(\bar{X},\bar{Y}) \\
 &=& -\eps  \G( X\wedge N + AX \wedge \phi, \J' \J  (Y\wedge N + AY \wedge \phi))\\
 &=& -\eps  \G( X\wedge N + AX \wedge \phi, \J'   (- \eps Y\wedge \phi + AY \wedge N))\\
 &=& -\eps  \G( X\wedge N + AX \wedge \phi,  - \eps ({\rm J'}Y) \wedge \phi + ({\rm J'}AY) \wedge N)\\
& =& -\eps  \left(\<X , {\rm J'}AY\>_p |N|_p^2 - \eps \<AX,{\rm J'}Y\>_p |\phi|^2_p \right)\\
&=  &- g(X, {\rm J}'AY)  + g(AX,{\rm J}'Y) \\
&=&g(X,(-{\rm J}'A+A {\rm J}')Y),
\end{eqnarray*}
which gives the claimed formula for $\bar{g}'.$ We now discuss the degeneracy and the signature of $\bar{g}',$ which depend on 
the type of the shape operator $A$:
\bigskip

\begin{itemize}
\item[-] \textbf{The real  diagonalizable case.} 
Write $g$ and $A$ in canonical form, with $(e_1,e_2)$ an oriented, orthonormal local frame. It follows that 
${\rm J}' e_1=e_2, {\rm J}'e_2 = -\eps' e_1$. We easily get
$$\bar{g}'= 
\left( \begin{array}{cc}
  0 & \eps_2 (\ka_2 -\ka_1) \\ 
 \eps_2 (\ka_2 -\ka_1) & 0
\end{array} \right), $$
We see in particular that $\bar{g}'$ is degenerate at umbilic points and indefinite otherwise.

\bigskip

\item[-] \textbf{The complex diagonalizable case.}
Write $g$ and $A$ in canonical form. It follows that
${\rm J}' e_1=e_2, {\rm J}'e_2 = e_1$ (here $\eps'=-1$ since the metric $g$ is indefinite). Hence
 $$\bar{g}'= 
\left( \begin{array}{cc}
 -2 \la& 0 \\ 
 0 & -2 \la
\end{array} \right), $$
 which shows that the metric $\bar{g}'$ is everywhere definite.

\item[-]  \textbf{The non diagonalizable case.}
Write $g$ and $A$ in canonical form. Since $(e_1,e_2)$ is a $g$-null basis, the complex structure is given 
by ${\rm J}' e_1=e_1, {\rm J}'e_2 = -e_2,$ so we get
 $$\bar{g}'= 
\left( \begin{array}{cc}
 0& 0 \\ 
 0 & -2
\end{array} \right), $$ which shows that the metric $\bar{g}'$ is everywhere degenerate. In particular, we don't need to take into consideration the case of $A$ being non diagonalizable in the proofs of Corollaries \ref{coroG'mini} and \ref{coroG'flat}.

\end{itemize}

\subsubsection{The mean curvature vector and the proof of Corollary \ref{coroG'mini}} \label{H'}
\noindent \textbf{The real diagonalizable case}

\noindent It has been seen in Section \ref{HG} that
$$\bar{h}_{ijj}:= \bar{h}(e_i,e_j,e_j)= -\eps \eps_j e_i(\kappa_j).$$
It follows that
$$ \G' (2\vec{H}',\J' d\bar{\phi}(e_1) )= \frac{\bar{h}_{112}}{\bar{g}'(e_1,e_2)}=
\frac{-\eps \eps_1 e_2(\ka_1)}{\eps_2(\ka_2 - \ka_1)}=-\eps \eps' \frac{ e_2(\ka_1)}{\ka_2 - \ka_1}$$
and
$$ \G' (2\vec{H}',\J' d\bar{\phi}(e_2) )= \frac{\bar{h}_{122}}{\bar{g}'(e_1,e_2)}
=\frac{-\eps \eps_2 e_1(\ka_2)}{\eps_2(\ka_2 - \ka_1)}=-\eps\frac{ e_1(\ka_2)}{\ka_2 - \ka_1}.$$
Hence
$$\vec{H}'= \frac{-\eps }{2(\ka_2-\ka_1)^2} \Big(\eps_1 e_1(\ka_2)  \J' d\bar{\phi}(e_1)+\eps_2 e_2(\ka_1) \J'  d\bar{\phi} (e_2) \Big). $$
In particular, we see that if $\bar{\s}$ is $\G'$-minimal, both $e_1(\ka_2)$ and $e_2(\ka_1)$ vanish. We now use the Codazzi equation derived in Section \ref{HG}:
$$ 
\left\{ \begin{array}{ccl}
 (\ka_2 -\ka_1) \omega_{12}^1 & = & \eps_2 e_1(\ka_2) \\ 
 (\ka_2 -\ka_1) \omega_{12}^2 & = & \eps_1 e_2(\ka_1).
\end{array} \right. $$ 
Since we assume that the metric $\bar{g}'$ is not degenerate, $\ka_2 -\ka_1$ does not vanish. Therefore  the $\G'$-minimality condition implies the
vanishing of $\omega_{12}^1$ and $\omega_{12}^2,$ i.e. the flatness of $g.$  The next step consists of using Gauss equation with respect to the immersion
$\phi: \M^2 \to \X^3_{p,1},$ giving

\begin{eqnarray*}
 g(R^g(e_1,e_2)e_1,e_2)&=& \eps  h(e_1,e_1)h(e_2,e_2)- \eps h(e_1,e_2)h(e_1,e_2) + K_{\X^3_{p,1}}\\
 &=&\eps \eps_1 \eps_2 \ka_1 \ka_2 + 1\\
 &=& \eps \eps' K +1.
 \end{eqnarray*}
Hence $ \ka_1 \ka_2 =-\eps \eps'.$ Taking into account  the vanishing of $e_1(\ka_2)$ and $e_2(\ka_1)$, it implies that both principal curvatures are constant, 
non vanishing and different of $\pm 1.$ In particular $\s$ has parallel second fundamental form and $\bar{\s}$ is totally geodesic.

\medskip

In the real diagonalizable case we are able to give a  more precise characterization of surfaces with parallel second fundamental form:
introducing  the map 
${\phi}_t:=\co (t) \phi + \si(t) N$ and differentiating, we
get
$$d{\phi}_t (e_2)= \co(t) d\phi(e_2) + \si(t) dNe_2=(\co (t)- \ka_2 \si (t)) d\phi(e_2).$$
Hence, choosing $t_0$ such that $\frac{\co (t_0)}{\si (t_0)}=\ka_2=-\eps(\ka_1)^{-1}$ yields the vanishing of $d{\phi}_{t_0} (e_2).$ Defining local coordinates $(s_1,s_2)$ on $\M^2$ such that $\pa_{s_1}=e_1$ and $\pa_{s_2}= e_2$, we claim that the curve $\ga(s_1):={\phi}_{t_0}(s_1,s_2)$  is a geodesic of $\X_{p,1}^3.$ 
To see this, we  calculate the  acceleration of $\ga$ in $\R^4$: 
\begin{eqnarray*}
\ga''&=&\frac{\pa^2 {\phi}_{t_0} }{\pa {s}_1^2} \\
& =&\co (t_0)\frac{\pa^2 \phi }{\pa s^2_1} + \si (t_0) \frac{\pa^2 N }{\pa s^2_1} \\
&=&\Big(\co (t_0) -\ka_1 \si (t_0) \Big)\frac{\pa^2 {\phi} }{\pa s^2_1}\\
&=& (\co (t_0) -\ka_1 \si (t_0) ) (\eps \eps_1 \ka_1 N - \eps_1 \phi)\\
&=&\frac{\co (t_0) -\ka_1 \si (t_0) }{\eps \co (t_0)}(-\si (t_0)N - \co(t_0) \phi),
\end{eqnarray*}
which is collinear to $\ga.$ Hence $\ga$ is a geodesic and $\phi(\M^2)$ is a tube over $\ga.$

\bigskip

\noindent \textbf{The complex diagonalizable case}

\noindent 
 Since the basis $(e_1,e_2)$ is orthogonal with respect to $\bar{g}'$, 
 the $\G'$-minimality of $\bar{\s}$ is equivalent to the vanishing of
 $$\bar{h}_{111}+ \bar{h}_{122}=-4 \eps \la \omega_{12}^1$$ 
 and 
   $$\bar{h}_{112}+ \bar{h}_{222}=-4 \eps \la \omega_{12}^2$$
(the coefficients $\bar{h}_{ijk}$ have been determined in Section \ref{HG2}).
Hence $\omega_{12}^1$ and  $\omega_{12}^2$ vanish and $g$ is flat. 
Again we use Gauss equation with respect to the immersion
$\phi: \M^2 \to \X^3_{p,1},$  obtaining
$$ g(R^g(e_1,e_2)e_1,e_2)= \eps h(e_1,e_1)h(e_2,e_2)- \eps h(e_1,e_2)h(e_1,e_2) + K_{\X^3_{p,1}}=-\eps(H^2+\la^2) +1$$
hence $ H^2+\la^2 =\eps.$   On the other hand, the Codazzi equation becomes a Cauchy-Riemann system satisfied by the pair $(H,\la)$:
$$\left\{  \begin{array}{ccr} e_1(H)&=&-e_2(\la)\\ e_2(H)&=&e_1(\la), \end{array}\right.$$
     so by Liouville theorem, $H$ and $\la$ are constant,
 which implies the vanishing of $\bar{h}.$
 
\subsubsection{Flat Lagrangian surfaces: proof of Corollary \ref{coroG'flat}}
\noindent \textbf{The real diagonalizable case}

 \noindent In order to characterize the flatness
 of $\bar{g}':=\bar{\phi}^* \G',$ we shall use twice the Gauss equation, first with respect to the immersion $\bar{\phi}: \M^2 \to L^{\pm}(\X^3_{p,1})$, and then with respect to the embedding $\iota:  L^{\pm}(\X^3_{p,1}) \to \Lambda^2(\R^4).$
 
 First, using the principal frame $(e_1,e_2)$ introduced in the previous section, we have
  $$K^{\bar{g}'} = \bar{g}'(R^{\bar{g}'}(\bar{e}_1,\bar{e}_2)\bar{e}_1,\bar{e}_2) = 
\G'(\vec{h}(\bar{e}_1,\bar{e}_2),\vec{h}(\bar{e}_1,\bar{e}_2) ) -\G'(\vec{h}(\bar{e}_1,\bar{e}_1),\vec{h}(\bar{e}_2,\bar{e}_2) )+\G'(\bar{R}(\bar{e}_1,\bar{e}_2)\bar{e}_1,\bar{e}_2),$$
where $\vec{h} : T \bar{\s} \times T \bar{\s} \to N \bar{\s}$ denotes the second fundamental form of the immersion $\bar{\phi}$ with respect to the metric $\G'.$
In other words,
$\G'(\vec{h}(X,Y),\J Z)=\bar{h}(X,Y,Z).$
We have
$$ \vec{h}(\bar{e}_i,\bar{e}_j)=\frac{\bar{h}_{ij2}N_1+\bar{h}_{ij1}N_2}{\eps_1 (\ka_2-\ka_1)},$$
so that
\begin{eqnarray*}\G'(\vec{h}(\bar{e}_1,\bar{e}_2),\vec{h}(\bar{e}_1,\bar{e}_2) )&=&2\eps_1 \frac{\bar{h}_{112} \bar{h}_{122}}{\ka_2-\ka_1}
\end{eqnarray*}
and
\begin{eqnarray*}\G'(\vec{h}(\bar{e}_1,\bar{e}_1),\vec{h}(\bar{e}_2,\bar{e}_2) )&=&\eps_1 \frac{\bar{h}_{111} \bar{h}_{222}+\bar{h}_{112} \bar{h}_{122}}{\ka_1-\ka_2}.
\end{eqnarray*}
Hence
\begin{eqnarray*}\G'(\vec{h}(\bar{e}_1,\bar{e}_2),\vec{h}(\bar{e}_1,\bar{e}_2) )-\G'(\vec{h}(\bar{e}_1,\bar{e}_1),\vec{h}(\bar{e}_2,\bar{e}_2) )
 &=&\eps_1 \frac{2\bar{h}_{112} \bar{h}_{122}-\bar{h}_{111} \bar{h}_{222}+\bar{h}_{112} \bar{h}_{122}}{\ka_1-\ka_2}\\
 &=&\eps_2 \frac{ e_2(\ka_1)e_1(\ka_2)- e_1(\ka_1)e_2(\ka_2)}{\ka_1-\ka_2}\\
 &=&-\eps_2 \frac{  (d\ka_1 \wedge d\ka_2) (e_1,e_2) }{\ka_1-\ka_2}.
\end{eqnarray*}
We now proceed to calculate $\G'(\bar{R}(\bar{e}_1,\bar{e}_2)\bar{e}_1,\bar{e}_2).$ We have
$$\bar{e}_i=d\phi(e_i) \wedge N + \phi \wedge dN(e_i)= -E_{2+i}-\ka_i E_i.$$
Then we easily get that $  h^{\iota}(\bar{e}_1,\bar{e}_1)= -\eps_1 (\eps +\ka_1^2)\bar{x}$ and $  h^{\iota}(\bar{e}_1,\bar{e}_2)=(\ka_1 - \ka_2 ) e_1 \wedge e_2.$ Analgously, $ h^{\iota}(\bar{e}_1, \eps \J' \circ \J \bar{e}_2) $ is collinear to $\bar{x}$, while $  h^{\iota}(\bar{e}_2, \eps \J' \circ \J \bar{e}_2)$
is collinear to $e_1 \wedge e_2.$

It follows that, using again Gauss equation and  the fact that the metric $\<\<.,.\>\>$ is flat,
\begin{eqnarray*}
 \G'(\bar{R}(\bar{e}_1,\bar{e}_2)\bar{e}_1,\bar{e}_2)&=& -\G'(\bar{R}(\bar{e}_1,\bar{e}_2)\bar{e}_1,\eps \J' \circ \J \bar{e}_2)\\
 &=&- \Big(\<\<h^{\iota}(\bar{e}_1, \J' \circ \J \bar{e}_2, h^{\iota}(\bar{e}_2,\bar{e}_1))\>\> 
 - \<\<h^{\iota}(\bar{e}_1,\bar{e}_1),h^{\iota}(\bar{e}_2, \eps \J' \circ \J' \bar{e}_2)\>\> \Big) \\
 &=& 0.
\end{eqnarray*}
We conclude that the metric $\bar{g}'$ (and therefore $\bar{g}$ as well) is flat if and only if $d\ka_1 \wedge d\ka_2$ vanishes, i.e.\
$\s$ is Weingarten.

\bigskip

\noindent \textbf{The complex diagonalizable case}

\noindent Since $\G'(N_i,N_i)=-\G'(\bar{e}_i,\bar{e}_i)=2\la,$ we have
$$ \vec{h}(\bar{e}_i,\bar{e}_j)=\frac{\bar{h}_{ij1}N_1+\bar{h}_{ij2}N_2}{2\la}.$$
Hence
\begin{eqnarray*}\G'(\vec{h}(\bar{e}_1,\bar{e}_2),\vec{h}(\bar{e}_1,\bar{e}_2))
 -\G'(\vec{h}(\bar{e}_1,\bar{e}_1),\vec{h}(\bar{e}_2,\bar{e}_2) )
 &=& -\frac{ \bar{h}_{112}^2+ \bar{h}_{122}^2- \bar{h}_{111} \bar{h}_{122}- \bar{h}_{112} \bar{h}_{222}}{2 \la}\\
 &=&  -\frac{  \bar{h}_{112}( \bar{h}_{112}- \bar{h}_{222})+  \bar{h}_{112}(\bar{h}_{112}- \bar{h}_{111})}{2 \la}\\
 &=&  -\frac{ 2e_1(\la)e_2(H)-2e_2(\la)e_1(H)}{2\la} \\
 &=&  \frac{ (dH \wedge d\la) (e_1,e_2)}{\la}
\end{eqnarray*}
It remains to prove that $\G'(\bar{R}(\bar{e}_1,\bar{e}_2)\bar{e}_1,\bar{e}_2)$ vanishes. We have
$$\bar{e}_1=d\phi(e_1) \wedge N + \phi \wedge dN(e_1)= -E_{3}-H E_1+ \la E_2$$
and
$$\bar{e}_2=d\phi(e_2) \wedge N + \phi \wedge dN(e_2)= -E_{4}-\la E_1 - H E_2 ,$$
so
$$\eps \J' \circ \J \bar{e}_2 = E_1 - \eps \la E_4 -\eps H E_3.$$
Then we easily get
\begin{eqnarray*}  h^{\iota}(\bar{e}_1,\bar{e}_1)&=& -\eps_1(1  +\eps (H^2 - \la^2)\bar{x}+ 2\la e_1 \wedge e_2\\
  h^{\iota}(\bar{e}_1,\eps \J' \circ \J' \bar{e}_2)&=& 0 \end{eqnarray*}
and
  $$  h^{\iota}(\bar{e}_2,\eps \J' \circ \J' \bar{e}_2)=(-1  -\eps (H^2 - \la^2) )e_1 \wedge e_2 + 2\la \eps \eps_1 \bar{x}.$$
 It follows that
\begin{eqnarray*}
 \G'(\bar{R}(\bar{e}_1,\bar{e}_2)\bar{e}_1,\bar{e}_2)&=& -\G'(\bar{R}(\bar{e}_1,\bar{e}_2)\bar{e}_1,\eps \J' \circ \J \bar{e}_2)\\
 &=&- \Big(\<\<h^{\iota}(\bar{e}_1,\eps \J' \circ \J \bar{e}_2), h^{\iota}(\bar{e}_2,\bar{e}_1)\>\> 
 - \<\<h^{\iota}(\bar{e}_1,\bar{e}_1),h^{\iota}(\bar{e}_2, \eps \J' \circ \J' \bar{e}_2)\>\> \Big) \\
 &=& 0.
\end{eqnarray*}
 We conclude that the metric $\bar{g}'$ (and therefore $\bar{g}$ as well) is flat if and only if $dH \wedge d\la$ vanishes, i.e.\
$\s$ is Weingarten.
 
\subsection{Marginally trapped Lagrangian surfaces: proof of Corollary \ref{coromargtrapped}}
\subsubsection{$\G$-marginally trapped Lagrangian surfaces}

We have seen in Section \ref{nondiagoG} that if the shape operator $A$ of $\phi$ is not diagonalizable, then $\bar{g}(e_1,e_1)$ vanishes. If follows that
$d\bar{\phi}(e_1)$, and therefore $\J d\bar{\phi}(e_1)$ as well, is a $\G$-null vector. We have also seen that
$\G(2 \vec{H},\J d\bar{\phi}(e_1))$ vanishes, so $\vec{H},$ a vector of the plane $N\bar{\s}$ spanned by $\J d\bar{\phi}(e_1)$ and $\J d\bar{\phi}(e_2)$, must be collinear to $\J d\bar{\phi}(e_1)$. Hence it is a $\G$-null vector as well. 
 
\subsubsection{$\G'$-marginally trapped Lagrangian surfaces}

We start from the expression of the mean curvature vector of $\bar{\s}$ with respect to $\G'$ obtained in Section \ref{H'}:
$$\vec{H}'= \frac{-\eps }{2(\ka_2-\ka_1)^2} \Big(\eps_1 e_1(\ka_2)  \J' d\bar{\phi}(e_1)+\eps_2 e_2(\ka_1) \J'  d\bar{\phi} (e_2) \Big) .$$
Since $\bar{g}'(e_1,e_1)$ and $\bar{g}'(e_1,e_1)$ vanish, the pair $(\J' d\bar{\phi}(e_1),\J' d\bar{\phi}(e_2))$ is a $\G$-null basis of the normal space $N \bar{\s}.$ Therefore, the mean curvature vector $\vec{H}'$ is $\G'$-null
 if and only if it is collinear to one of the two vectors $\J' d\bar{\phi}(e_i)$, i.e.\ if and only if either
$e_1(\ka_2)$ or $e_2(\ka_1)$ vanishes. This occurs at least in the following two cases:
\begin{itemize}
	\item[-] If $\s$ is a tube, i.e.\ the set of equidistant points to a given curve of $\X^3_{p,1}$, then one of its principal curvatures is constant;
	\item[-] If $\s$ is a surface of revolution, i.e. a surface invariant by the action of a subgroup $SO(2)$ or $SO(1,1)$ of $SO(4-p,p)$, then both principal curvatures are constant along the orbits of the action, which are in addition tangent to one of the principal directions (cf \cite{An}). Therefore, $e_1(\ka_2)$ or $e_2(\ka_1)$ vanishes.
\end{itemize}

\end{document}